\documentclass{article}
\usepackage{hyperref}
\usepackage{array,float,bigstrut,amsmath,amsthm,amssymb,mathrsfs}
\usepackage[plain]{algorithm}
\usepackage{algpseudocode}

\newcommand{\A}{\mathscr{A}}
\newcommand{\B}{\mathscr{B}}
\newcommand{\E}{\mathcal{E}}
\newcommand{\G}{\mathscr{G}}
\newcommand{\Q}{\mathcal{Q}}
\renewcommand{\P}{\mathcal{P}}
\newcommand{\R}{\mathcal{R}}
\renewcommand{\S}{\mathcal{S}}
\newcommand{\QP}{(\Q,\P)}
\newcommand{\K}{\mathcal{K}}
\newcommand{\T}{\mathcal{T}}
\newcommand{\SN}{\mathscr{N}}
\newcommand{\SP}{\mathscr{P}}

\renewcommand{\>}{\rangle}
\newcommand{\cl}{\mathrm{cl}}
\newcommand{\opts}{\mathrm{opts}}
\newcommand{\mex}{\mathrm{mex}}

\newtheorem{theorem}{Theorem}[section]
\newtheorem{definition}[theorem]{Definition}

\newtheorem{corollary}[theorem]{Corollary}
\newtheorem{lemma}[theorem]{Lemma}
\newtheorem{fact}[theorem]{Fact}

\newcommand{\PPos}{$\SP$-position}
\newcommand{\NPos}{$\SN$-position}
\newcommand{\Pow}{\textit{Pow}}

\newcommand{\vsp}{\vspace{0.15cm}}

\newcommand{\Star}{\ast}

\input macros.tex

\begin{document}
\title{The Structure and Classification\\of Mis\`ere Quotients\\ {\small PREPRINT}}
\author{Aaron N. Siegel\\Institute for Advanced Study\\1 Einstein Drive\\Princeton, NJ 08540}
\date{\today}
\maketitle

\begin{abstract}
A \emph{bipartite monoid} is a commutative monoid $\Q$ together with an identified subset $\P \subset \Q$.  In this paper we study a class of bipartite monoids, known as \emph{mis\`ere quotients}, that are naturally associated to impartial combinatorial games.

We introduce a structure theory for mis\`ere quotients with $|\P| = 2$, and give a complete classification of all such quotients up to isomorphism.  One consequence is that if $|\P| = 2$ and $\Q$ is finite, then $|\Q| = 2^n+2$ or $2^n+4$.

We then develop computational techniques for enumerating mis\`ere quotients of small order, and apply them to count the number of non-isomorphic quotients of order at most~$18$.  We also include a manual proof that there is exactly one quotient of order~$8$.
\end{abstract}

\section{Introduction}
\label{section:introduction}

An \emph{impartial combinatorial game} $\Gamma$ is a two-player game with no hidden information and no chance elements, in which both players have exactly the same moves available at all times.  When $\Gamma$ is played under the \emph{mis\`ere-play condition}, the player who makes the last move loses.

Thirty years ago, Conway~\cite{conway_1976} showed that the mis\`ere-play combinatorics of such games are often frighteningly complicated.  However, new techniques recently pioneered by Plambeck~\cite{plambeck_2005} have reinvigorated the subject.  At the core of these techniques is the \emph{mis\`ere quotient}, a commutative monoid that encodes the additive structure of an impartial combinatorial game (or a set of such games).  See~\cite{siegel_mqlectures} for a gentle introduction to mis\`ere quotients, and~\cite{siegel_200Xd} for a more rigorous one; see~\cite{plambeck_200X} for a survey of the theory.

The introduction of mis\`ere quotients opens up a fascinating new area of study: the investigation of their algebraic properties.  Such investigations are intrinsically interesting, and also have the potential to reveal new insights into the mis\`ere-play structure of combinatorial games.  In this paper, we introduce several new results that expose quite a bit of structure in mis\`ere quotients.

Henceforth, we assume familiarity with mis\`ere quotients, and in particular with the recent work of Plambeck and Siegel~\cite{siegel_200Xd}.

\subsection*{Tame Extensions}

The first and most striking result concerns mis\`ere quotients with $\SP$-portions of size~2.  If $\QP$ is a mis\`ere quotient, then the \emph{tame extension} $\T\QP$ is a certain conservative extension of $\QP$ that adds no new $\SP$-positions.  It is defined in such a way that
\[\T_3 = \T(\T_2),\ \T_4 = \T(\T_3),\ \T_5 = \T(\T_4),\ \ldots\]
If we replace the ``base'' $\T_2$ by another quotient, say $\R_8$, we get a new family
\[\R_8,\ \T(\R_8),\ \T(\T(\R_8)),\ \ldots\]
and since $\R_8$ has a size-2 $\SP$-portion, so does every quotient in the new sequence.  The main result is that \emph{every} finite quotient with $|\P| = 2$ is isomorphic to a quotient in one of these two families.  It will follow that every finite quotient with $|\P| = 2$ has order $2^n+2$ or $2^n+4$, for some $n \geq 2$. Furthermore, if $\QP$ is an \emph{infinite} quotient with $|\P| = 2$, then $\QP \cong$ either $\T_\infty$ or $\R_\infty$, the limits of the two families (in a sense to be precisely defined).

\subsection*{``Almost Tame'' Octal Games}

Tame extensions also have a useful (and quite beautiful) application to octal games.  Fix an octal game~$\Gamma$ and an integer~$M$, and consider the partial quotient $\QP = \Q_M(\Gamma)$ and pretending function~$\Phi_M$.  Assume that $\QP$ is normal and~$\Q$ is finite, and let $\K \subset \Q$ be the kernel of $\Q$ (i.e., the intersection of all ideals in $\Q$).  In Section~\ref{section:normality}, we will show that if $\Phi_M(H_n) \in \K$ for sufficiently many heaps~$H_n$, then $\Q_{M+1}(\Gamma)$ is either $\QP$ or $\T\QP$, and $\Phi_{M+1}(H_{M+1}) \in \K$.  ``Sufficiently many'' will be in the Guy--Smith sense.

This theorem can be iterated, with strong consequences.  In particular, if we determine that $\Phi_M(H_n) \in \K$ for sufficiently many~$n$, then we can conclude that $\Q(\Gamma)$ is one of
\[\QP,\ \T\QP,\ \T(\T\QP),\ \ldots,\]
or possibly the limit $\T^\infty\QP$ of this sequence.  Furthermore, $\Phi(H_n) \in \K$ for all but finitely many~$n$.  Since $\Q(\Gamma)$ is normal, normal and mis\`ere play coincide on~$\K$; so we conclude that mis\`ere play reduces to normal play unless all the heaps are small.  In practice, this means that once we have computed $\Q_M(\Gamma)$, then we have completely characterized the ``mis\`ere-play divergence'' of~$\Gamma$; and its mis\`ere-play solution now depends only on finding a normal-play solution.

An ideal example is the game \textbf{0.414}, which we mentioned in~\cite{siegel_200Xd}.  Its normal-play solution is unknown, despite the computation of at least $2^{24}$ $\G$-values by Flammenkamp~\cite{flammenkamp_www_octal}.  However, it is easy enough to compute $\Q_{18}(\textbf{0.414})$, and to verify using the above logic that $\Phi(H_n) \in \K$ for all $n > 18$.  Thus we know $\Q(\textbf{0.414}) \cong \T^k(\Q_{18})$, for some~$k$, and we need invest no further worry in the mis\`ere play of \textbf{0.414}: we may sit back and await a normal-play solution.

One might recall the mis\`ere-play strategy for \textsc{Nim}: Play normal \textsc{Nim} unless your move would leave only heaps of size~$1$.  In that case, play to leave an odd number of heaps of size~$1$.  We can now state an analogous strategy for \textbf{0.414}: Play normal \textbf{0.414} unless your move would leave only heaps of size $\leq 18$.  In that case, consult the fine structure of~$\Q_{18}$.  We can state this reduction with confidence, despite the fact that the normal-play strategy for \textbf{0.414} remains unknown.

Section~\ref{section:limits} is mostly spadework.  In Section~\ref{section:normality}, we define tame extensions, prove a key result showing that certain extensions are always tame, and apply this result to octal games.  In Sections~\ref{section:Tn} and~\ref{section:Rn}, we develop a structure theory for the quotients $\T_n$ and $\T^n(\R_8)$, and use this machinery to prove the main theorem on $|\P| = 2$ quotients.

\subsection*{Quotients of Small Order}
\suppressfloats[t]

We prove in Section~\ref{section:R8} that $\R_8$ is the only quotient of order~8 (up to isomorphism), but the primary effect of that proof is to discourage attempts to extend this classification by hand.  The rest of the paper focuses instead on developing computational techniques for classifying quotients of small order.

In Section~\ref{section:validity}, we show that an arbitrary r.b.m.\ $\QP$ is a mis\`ere quotient if and only if there exists a \emph{valid transition table} for $\QP$---a certain combinatorial structure superimposed on $\QP$.  This yields a computational method for testing whether $\QP$ is a mis\`ere quotient, which is optimized and applied in Section~\ref{section:enumeration}.  The fruits of this effort are summarized in Figure~\ref{figure:smallorder}.

\begin{figure}[htbp]
\centering
\begin{tabular}{r|rrrrrrrrrr}
$n$                       & 2 & 4 & 6 & 8 & 10 & 12 & 14 & 16 &  18 \bigstrut \\ \hline
Quotients of order $n$    & 1 & 0 & 1 & 1 &  1 &  6 &  9 & 50 & 211 \bigstrut
\end{tabular}
\caption{The number of mis\`ere quotients of order $n \leq 18$ (up to isomorphism).}
\label{figure:smallorder}
\end{figure}

\subsection*{Preliminaries}

We recall some key facts and definitions from~\cite{siegel_200Xd}, and also introduce some new notation.

Let $\QP$ be a bipartite monoid.  Two elements $x,y \in \Q$ are \emph{indistinguishable} if, for all $z \in \Q$,
$xz \in \P \Leftrightarrow yz \in \P$.  $\QP$ is \emph{reduced} if the elements of $\Q$ are pairwise distinguishable.  In~\cite{siegel_200Xd} we showed that every bipartite monoid has a unique reduced quotient.

We say that $\QP$ is a \emph{sub-b.m.} of $(\S,\R)$ if $\Q$ is a submonoid of $\S$ and $\R \cap \Q = \P$.  In this case we write $\QP < (\S,\R)$.


\begin{definition}
Let $\QP$ be a bipartite monoid and fix $x \in \Q$.  The \emph{meximal set of $x$ in $\QP$}, denoted $\mathcal{M}_x$, is defined by
\[\mathcal{M}_x = \{y \in \Q : \textrm{there is no } z \in \Q \textrm{ such that } xz,yz \in \P\}.\]
\end{definition}

The following statement is slightly more general than the rule given in~\cite{siegel_200Xd}, but the proof is identical.

\begin{fact}[Generalized Mex Rule]
\label{fact:generalizedmexrule}
Let $\QP = \Q(\A)$, and let $\QP < (\S,\R)$.  Fix $G$ with $\opts(G) \subset \A$ and fix $x \in \S$.  The following are equivalent.
\begin{enumerate}
\item[(a)] $\Q(\A \cup \{G\}) \cong (\S,\R)$ and $\Phi(G) = x$.
\item[(b)] $\S$ is generated by $\Q \cup \{x\}$, and the following two conditions hold.
\begin{enumerate}
\item[(i)] $\Phi''G \subset \mathcal{M}_x$; and
\item[(ii)] For each $Y \in \A$ and $n \geq 0$ such that $x^{n+1}\Phi(Y) \not\in \P$, we have either: $x^{n+1}\Phi(Y') \in \P$ for some option $Y'$ of $Y$; or else $x^nx'\Phi(Y) \in \P$ for some $x' \in \Phi''G$.
\end{enumerate}
\end{enumerate}
\end{fact}

In~\cite{siegel_200Xd} we stated Fact~\ref{fact:generalizedmexrule} for the special case $\S = \Q$.  This will often be the case of greatest interest, but we shall have several occasions to use the more general form.

\section{Limits and One-Stage Extensions}
\label{section:limits}

In this section we show that every mis\`ere quotient is the limit of a sequence of finitely generated quotients.  Furthermore, each term of this sequence is a conservative extension of the previous term, in a way we now make precise.

\begin{definition}
Let $\QP$, $(\Q^+,\P^+)$ be reduced bipartite monoids.  We say that $(\Q^+,\P^+)$ is an \emph{extension} of $\QP$ if there is some submonoid $(\S,\R) < (\Q^+,\P^+)$ such that $\QP$ is (isomorphic to) the reduction of $(\S,\R)$.  If $\Q^+$ is generated by $\S \cup \{x\}$ for some single element $x \in \Q^+ \setminus \S$, then we say that $(\Q^+,\P^+)$ is a \emph{one-stage extension} of $\QP$.
\end{definition}

\begin{lemma}
\label{lemma:stages}
Let $\QP$ be a finitely generated mis\`ere quotient.  Then there is a sequence of mis\`ere quotients
\[0 = (\Q_0,\P_0),\ (\Q_1,\P_1),\ \ldots,\ (\Q_n,\P_n) = \QP\]
such that each $(\Q_{i+1},\P_{i+1})$ is a one-stage extension of $(\Q_i,\P_i)$.
\end{lemma}

\begin{proof}
Write $\QP = \Q(\A)$ and choose a finite set $\mathscr{H} \subset \A$ so that $\Phi''\mathscr{H}$ generates $\Q$.  Since the hereditary closure of a finite set is finite, we may assume that $\mathscr{H}$ is hereditarily closed.  Enumerate
\[\mathscr{H} = \{H_0,H_1,\ldots,H_m\}\]
so that the successive $H_i$'s have nondecreasing birthdays, and put
\[(\Q_i,\P_i) = \Q(H_0,\ldots,H_i).\]
It is easily seen that either $(\Q_{i+1},\P_{i+1}) = (\Q_i,\P_i)$, or else it is a one-stage extension of $(\Q_i,\P_i)$.  A suitable reindexing gives the lemma.
\end{proof}

Now let $(\Q_n,\P_n)$ be a sequence of bipartite monoids, and for each $n$, let $(\Q_n^+,\P_n^+) < (\Q_{n+1},\P_{n+1})$ and let $\pi_n : \Q_n^+ \to \Q_n$ be a surjective homomorphism of bipartite monoids.  We call $(\Q_n,\P_n,\pi_n)$ a \emph{partial inverse system}.

Let $\overleftarrow{\Q} = (\Q_n,\P_n,\pi_n)$ be a partial inverse system.  It is convenient to regard the underlying sets of the $\Q_n$ as formally disjoint.  A \emph{thread of $\overleftarrow{\Q}$ starting at $n$} is a sequence $(x_n,x_{n+1},x_{n+2},\ldots)$, where $x_n \in \Q_n$ and for each $i > n$ we have $x_{i+1} \in \Q^+_i$ and $\pi_i(x_{i+1}) = x_i$.  We say two threads $\overleftarrow{x}$ and $\overleftarrow{y}$ are equivalent, and write $\overleftarrow{x} \sim \overleftarrow{y}$, if one is a terminal segment of the other.

If $\overleftarrow{x} = (x_m,x_{m+1},x_{m+2},\ldots)$ and $\overleftarrow{y} = (y_n,y_{n+1},y_{n+2},\ldots)$ are threads, we can define their product as follows.  Without loss of generality, assume that $m \leq n$, and put
\[\overleftarrow{x} \cdot \overleftarrow{y} = (x_ny_n,x_{n+1}y_{n+1},x_{n+2}y_{n+2},\ldots)\]
It is easy to check that $\overleftarrow{x} \cdot \overleftarrow{y}$ is a thread and that the product respects the equivalence $\sim$.  Further, $\overleftarrow{1} \cdot \overleftarrow{x} = \overleftarrow{x}$, where $\overleftarrow{1} = (1,1,1,\ldots)$ is a list of the identity elements of each $\Q_n$.  Thus the threads modulo $\sim$ form a commutative monoid $\Q$.  We can define a subset $\P \subset \Q$ by
\[\P = \{(x_n,x_{n+1},x_{n+2},\ldots) \in \Q : \textrm{ some (all) } x_i \in \P_i\},\]
and this makes $(\Q,\P)$ into a bipartite monoid, which we call the \emph{partial inverse limit} of the system $\overleftarrow{Q}$.  We write $(\Q,\P) = \lim \overleftarrow{Q} = \lim_n(\Q_n,\P_n)$.

The following lemma is an easy exercise.

\begin{lemma}
\label{lemma:reducedlimit}
If $(\Q_n,\P_n)$ is reduced for infinitely many values of $n$, then so is $\lim_n(\Q_n,\P_n)$.
\end{lemma}

\begin{theorem}
\label{theorem:chainlimit}
Suppose that $\A_0 \subset \A_1 \subset \A_2 \subset \cdots$ is a chain of closed sets of games.  Then the quotients $\Q(\A_n)$ form a partial inverse system, and we have
\[\Q\left(\bigcup_n \A_n\right) \cong \lim_n\Q(\A_n).\]
\end{theorem}

\begin{proof}
Let $\Phi_n : \A_n \to \Q(\A_n)$ be the quotient maps and put $\Q_n^+ = \Phi_{n+1}''\A_n$.  Define $\pi_n : \Q_n^+ \to \Q_n$ by $\pi_n(\Phi_{n+1}(X)) = \Phi_n(X)$.  Now if $X \equiv_{\A_{n+1}} Y$, then necessarily $X \equiv_{\A_n} Y$, so $\pi_n$ is well-defined.

Now by Lemma~\ref{lemma:reducedlimit}, $\lim_n\Q(\A_n)$ is reduced.  To complete the proof, it suffices to exhibit a surjective homomorphism $\Phi : \bigcup_n\A_n \to \lim_n\Q(\A_n)$.  Let $n$ be least so that $X \in \A_n$, and put
\[\Phi(X) = (\Phi_n(X),\Phi_{n+1}(X),\Phi_{n+2}(X),\ldots).\]
It is easily verified that $\Phi$ has the desired properties.
\end{proof}

An easy corollary of Theorem \ref{theorem:chainlimit} will be central to the classification theory.

\begin{corollary}
\label{corollary:fglimit}
Suppose that $(\Q,\P)$ is a non-f.g.\ mis\`ere quotient.  Then there is some partial inverse system $(\Q_n,\P_n)$ of finitely generated mis\`ere quotients such that:
\begin{itemize}
\item[(i)] $(\Q_0,\P_0) = 0$;
\item[(ii)] Each $(\Q_{n+1},\P_{n+1})$ is a one-stage extension of $(\Q_n,\P_n)$; and
\item[(iii)] $(\Q,\P) = \lim_n(\Q_n,\P_n)$.
\end{itemize}
\end{corollary}

\begin{proof}
Write $(\Q,\P) = \Q(\A)$ with $\A$ closed.  Enumerate $\A = \{H_0,H_1,H_2,\ldots\}$ so that the birthdays of the $H_n$ are nondecreasing.  (This can always be done, since there are only finitely many games of each fixed birthday.)  Then for each~$n$, we have $\opts(H_n) \subset \{H_0,\ldots,H_{n-1}\}$.  Put
\[(\Q_n,\P_n) = \Q(H_0,\ldots,H_n).\]
Let $\Q_n^+$ be the submonoid of $\Q_{n+1}$ generated by $\{[H_0],\ldots,[H_n]\}$, and define $\pi_n : \Q_n^+ \to \Q_n$ by
\[\pi_n([H]_{\A_{n+1}}) = [H]_{\A_n}.\]
$\pi_n$ is well-defined, since each $G \equiv_{\A_{n+1}} G'$ implies $G \equiv_{\A_n} G'$.  Now (i) is immediate, since necessarily $H_0 = 0$, and (ii) follows easily (after reindexing to eliminate cases where $\Q_{n+1} = \Q_n$).  Now by Lemma~\ref{lemma:reducedlimit}, we know that $\lim_n(\Q_n,\P_n)$ is reduced.  To prove (iii), it therefore suffices to show that $\lim_n(\Q_n,\P_n)$ is a quotient of $\A$.

Let $\Phi_n : \cl(\{H_0,\ldots,H_n\}) \to \Q_n$ be the usual quotient map, and define $\Phi : \A \to \Q$ by
\[\Phi(H_n) = (\Phi_n(H_n),\Phi_{n+1}(H_n),\Phi_{n+2}(H_n),\ldots).\]
It is easily checked that $\Phi$ is a surjective homomorphism of bipartite monoids.
\end{proof}

\section{Normal Quotients and Tame Extensions}
\label{section:normality}

In this section we introduce a certain algebraic property known as \emph{faithful normality}, and we study one-stage extensions of faithfully normal quotients.  In particular, we show that certain one-stage extensions of faithfully normal quotients behave exactly like normal-play Grundy extensions.  The vast majority of quotients encountered in practice are faithfully normal, so this work has useful applications to octal games.

\begin{definition}
Let $\QP$ be a mis\`ere quotient with kernel $\K$, and let $z \in \K$ be the kernel identity.  We say that $\QP$ is \emph{regular} if $|\K \cap \P| = 1$, and \emph{normal} if $\K \cap \P = \{z\}$.
\end{definition}

\begin{definition}
Let $\QP = \Q(\A)$ and let $\Phi : \A \to \Q$ be the quotient map.  Suppose that
\[\Phi(G) = \Phi(H) \Longrightarrow \G(G) = \G(H) \qquad \textrm{for all } G,H \in \A.\]
Then we say that $\Phi$ is \emph{faithful}.  If in addition $\QP$ is normal, then we say that $\Phi$ is \emph{faithfully normal}.
\end{definition}

Often we will abuse terminology and refer to the \emph{quotient} as being faithful (or faithfully normal), rather than the quotient map.  We recall the following fact from~\cite{siegel_200Xd}.

\begin{fact}
Suppose $\Phi : \A \to \Q$ is faithfully normal.  Then $\K$ is isomorphic to the normal quotient of $\A$.
\end{fact}

Roughly speaking, therefore, faithful normality asserts that normal and mis\`ere play coincide on~$\K$. We have $\K \cong \mathbb{Z}_2^n$ for some $n$, and for each $i < 2^n$ there is a unique $z_i \in \K$ representing games of Grundy value $i$.  For convenience, when $\E \subset \K$, we write $z_m = \mex(\E)$ to mean $m = \mex\{i : z_i \in \E\}$.

Now fix a faithfully normal quotient $\Q(\A)$ with kernel~$\K$, and let $G \neq 0$ be a game such that $\opts(G) \subset \A$.   Then $\Q(\A \cup \{G\})$ is necessarily a one-stage extension of $\Q(\A)$.  For the remainder of this section, we will focus on the special case where $\Phi''G \subset \K$.  We will show that in this case, one-stage extensions behave \emph{exactly} like normal-play Grundy extensions.  In particular:

\begin{itemize}
\item Extensions by a \emph{proper} subset of the kernel are conservative and follow the mex rule.  Formally, if $\Phi''G \subsetneqq \K$, then $\Q(\A \cup \{G\}) \cong \Q(\A)$ and $\Phi(G) = \mex(\Phi''G)$.
\item Extensions by the \emph{entire} kernel cause the kernel to grow (from $\mathbb{Z}_2^n$ to $\mathbb{Z}_2^{n+1}$).  They behave like normal-play extensions whose Grundy values are new powers of~2.  Formally, if $\Phi''G = \K$, then $\Q(\A \cup \{G\}) \cong \T(\Q(\A))$, where $\T(\Q(\A))$ is a certain ``tame extension'' of $\Q(\A)$ that generalizes the extension $\mathbb{Z}_2^n < \mathbb{Z}_2^{n+1}$.
\end{itemize}

We begin with the $\Phi''G \subsetneqq \K$ case.

\begin{lemma}
\label{lemma:kernelextensions}
Suppose $\QP = \Q(\A)$ is faithfully normal with kernel $\K$.  Let $G$ be a game with $\opts(G) \subset \A$ and suppose $\Phi''G \subsetneqq \K$.  Then $\Q(\A \cup \{G\}) \cong \Q(\A)$ and $\Phi(G) = \mex(\Phi''G)$.
\end{lemma}

\begin{proof}
We verify conditions (i) and (ii) of the Generalized Mex Rule, with $(\S,\R) = \QP$ and $x = \mex(\Phi''G)$.  Note that $x = z_m$, where $m = \G(G)$.

For (i), normality implies that $\mathcal{K} \setminus \{z_m\} \subset \mathcal{M}_x$.  Since $\Phi''G \subset \K \setminus \{z_m\}$, this suffices.  For~(ii), fix $Y \in \A$ and $n \geq 0$, and suppose $x^{n+1}\Phi(Y) \not\in \P$.  If $n$ is odd, then since $x \in \K$ and the quotient is faithfully normal, we have $\G(Y) > 0$.  Thus $\G(Y') = 0$ for some $Y'$, whence $x^{n+1}\Phi(Y') \in \P$.  

Conversely, suppose that $n$ is even.  Then $\G(Y) \neq m$.  If $\G(Y) > m$, then $\G(Y') = m$ for some $Y'$, whence $x^{n+1}\Phi(Y') \in \P$.  Otherwise, let $i = \G(Y)$.  Since $\Phi''G \subset \K$ and $z_m = \mex(\Phi''G)$, we necessarily have $z_i \in \Phi''G$.  But $z_i\Phi(Y) = z$, so $x^nz_i\Phi(Y) = z \in \P$.
\end{proof}

\subsection*{Tame Extensions}

We now consider the case where $\Phi''G = \K$.  Let $\QP$ be a bipartite monoid with kernel $\K$, and define
\[\overline{\K} = \{\overline{x} : x \in \K\},\]
where each $\overline{x}$ is taken to be a formal symbol.

\begin{definition}
The \emph{first tame extension} $\T\QP = (\Q^+,\P^+)$ is defined as follows.  $\Q^+ = \Q \cup \overline{\K}$, $\P^+ = \P$, and multiplication is extended by:
\[x \cdot \overline{y} = \overline{xy} \ (x \in \Q,\ y \in \K); \qquad
  \overline{x} \cdot \overline{y} = xy \ (x,y \in \K).\]
The \emph{$n^\mathrm{th}$ tame extension} $\T^n\QP$ is defined by
\[\T^0\QP = \QP; \qquad \T^{n+1}\QP = \T(\T^n\QP).\]
Finally, we define
\[\T^\infty\QP = \lim_n \T^n\QP.\]
\end{definition}

Observe that the sequence of normal quotients
\[0,\ \mathbb{Z}_2,\ \mathbb{Z}_2^2,\ \mathbb{Z}_2^3,\ \ldots,\ \mathbb{Z}_2^\mathbb{N}\]
can be written
\[\T^0(0),\ \T^1(0),\ \T^2(0),\ \T^3(0),\ \ldots,\ \T^\infty(0)\]
while the sequence of tame mis\`ere quotients
\[\T_0,\ \T_1,\ \T_2,\ \T_3,\ \T_4,\ \ldots,\ \T_\infty\]
can be written
\[\T_0,\ \T_1,\ \T^0(\T_2),\ \T^1(\T_2),\ \T^2(\T_2),\ \ldots,\ \T^\infty(\T_2)\]
Thus the normal quotients can be viewed as a tame sequence with base 0, and the tame mis\`ere quotients can be viewed as a tame sequence with base $\T_2$.

If $\QP$ is a mis\`ere quotient, then so is $\T\QP$, as the following lemma establishes (cf.~Lemma~\ref{lemma:kernelextensions}).

\begin{lemma}
\label{lemma:tameextensions}
Suppose $\QP = \Q(\A)$ is faithfully normal with kernel $\K$.  Let $G$ be a game with $\opts(G) \subset \A$ and suppose $\Phi''G = \K$.  Then $\Q(\A \cup \{G\}) \cong \T(\Q(\A))$ and $\Phi(G) = \overline{z}$.
\end{lemma}

\begin{proof}
Identical to the proof of Lemma~\ref{lemma:kernelextensions}.
\end{proof}

\begin{corollary}
Suppose $\Q(\A)$ is faithfully normal with kernel $\K$.  Then for all $n \in \mathbb{N} \cup \{\infty\}$, $\T^n(\Q(\A))$ is a mis\`ere quotient.
\end{corollary}

\begin{proof}
Let $G_0 = 0$ and $\A_0 = \A$.  Recursively choose $G_{n+1}$ so that $\opts(G_{n+1}) \subset \A_n$ and $\Phi''G_{n+1} = \ker \Q(\A_n)$.  Put $\A_{n+1} = \cl(\A_n \cup \{G_{n+1}\})$.

By repeated application of Lemma~\ref{lemma:tameextensions}, we have $\Q(\A_n) \cong \T^n(\Q(\A))$, and Theorem~\ref{theorem:chainlimit} therefore gives $\Q(\bigcup_n \A_n) = \T^\infty(\Q(\A))$.
\end{proof}

\subsection*{The Quotients $\mathcal{R}_{2^n+4}$}

If we start with a different base $\QP$, we obtain another sequence of quotients $\T^n\QP$.  For example, if $\QP = \mathcal{R}_8$, then for all $n \geq 2$, $\T^{n-2}\QP$ is a quotient of order $2^n + 4$, which we denote by $\mathcal{R}_{2^n+4}$.  Likewise, we define $\mathcal{R}_\infty = \T^\infty(\mathcal{R}_8)$.  Since $|\P| = 2$, all the $\mathcal{R}_n$'s have $\SP$-portions of size~2.  A major goal of this paper is to prove the following theorem.

\begin{theorem}
\label{theorem:pportionsize2}
Suppose $\QP$ is a mis\`ere quotient with $|\P| = 2$.  Then either $\QP \cong \T_n$ or $\QP \cong \mathcal{R}_n$, for some $n \in \mathbb{N} \cup \{\infty\}$.
\end{theorem}

Thus if $\QP$ is a mis\`ere quotient with $|\P| = 2$, it follows that either $|\Q| = \infty$, or $|\Q| = 2^n + 2$ or $2^n + 4$ for some $n \geq 2$.  Furthermore, there is exactly one such quotient of each permissible finite order, and exactly two such infinite quotients.

\subsection*{``Almost Tame'' Octal Games}

Lemmas \ref{lemma:kernelextensions} and \ref{lemma:tameextensions} have useful implications for octal games, as summarized by the following theorem.

\begin{theorem}
\label{theorem:almosttameness}
Let $\Gamma$ be an octal game with last non-zero code digit $d$.  Fix~$n_0$, and suppose that $\QP = \Q_{2n_0+d-1}(\Gamma)$ is faithfully normal with kernel~$\K$.  Suppose furthermore that
\[\Phi(H_n) \in \K \textrm{ for all $n$ such that } n_0 \leq n < 2n_0+d.\]
Then:
\begin{itemize}
\item[(i)] $\Q(\Gamma)$ is a faithfully normal quotient;
\item[(ii)] $\Q(\Gamma) \cong \T^k\QP$ for some $k \in \mathbb{N} \cup \{\infty\}$;
\item[(iii)] $\Phi(H_n) \in \ker\Q(\Gamma)$, for all $n \geq n_0$.
\end{itemize}
\end{theorem}

\begin{proof}
We first show that (i)-(iii) hold for $\Q_n(\Gamma)$, for all $n$.  By hypothesis we may assume that $n \geq 2n_0 + d$.  Then a typical option of $H_n$ is a position $H_a + H_b$, with $a + b \geq 2n_0$.  Without loss of generality, we have $a \geq n_0$, so by induction $\Phi(H_a) \in \K$.  Thus $\Phi(H_a)x \in \K$ for all $x$, and in particular $\Phi(H_a + H_b) \in \K$.  This shows that $\Phi''H_n \subset \K$, and Lemmas~\ref{lemma:kernelextensions} and~\ref{lemma:tameextensions} immediately imply (i)-(iii).

If the partial quotients $\Q_n(\Gamma)$ eventually converge to some $\T^k\QP$, then $\Q(\Gamma) \cong \T^k\QP$.  Otherwise $\Q(\Gamma) \cong \T^\infty\QP$; and in either case (i)-(iii) are immediate.
\end{proof}

Thus when the hypotheses of Theorem~\ref{theorem:almosttameness} are satisfied, we know that beyond heap~$n_0$, the mis\`ere-play analysis of $\Gamma$ is no harder than its normal-play analysis.  It follows that we can stop computing partial quotients of $\Gamma$ and revert to the much easier task of calculating Grundy values.  We may say that $\Gamma$ is \emph{tame relative to heap $n_0$}.

The hypotheses of Theorem~\ref{theorem:almosttameness} may seem rather restrictive, but there are several three-digit octal games that satisfy them; for example, \textbf{0.414}, \textbf{0.776}, and \textbf{4.76}.  The mis\`ere-play solutions to these games now depend only on finding normal-play solutions, and we can regard them as ``relatively solved.''

The hexadecimal game \textbf{0.9092} is another interesting case.  It is known to be arithmeto-periodic in normal play.  Furthermore, in mis\`ere play we can show that it is tame relative to heap~12.  Now $\Q_{12}(\textbf{0.9092}) \cong \R_8$, so by Theorem~\ref{theorem:almosttameness} (suitably generalized to hexadecimal games) we have $\Q(\textbf{0.9092}) \cong \T^k(\R_8)$ for some $k$.  Since the $\G$-values of \textbf{0.9092} are unbounded, $k$ is necessarily $\infty$.  Therefore $\Q(\textbf{0.9092})$ is exactly $\R_\infty$.

\section{One-Stage Extensions of $\mathcal{T}_n$}
\label{section:Tn}

We next focus our attention on proving Theorem~\ref{theorem:pportionsize2}.  The crux of the proof is an analysis of one-stage extensions of $\mathcal{T}_n$ and $\mathcal{R}_{2^n+4}$.  This analysis also yields a useful structure theory for these quotients.  In particular, we will prove the following two theorems.

\begin{theorem}
\label{theorem:tnextensions}
If $\QP$ is a one-stage extension of $\T_n$ and $|\P| = 2$, then either $\QP \cong \T_{n+1}$, or else $\QP \cong \R_{2^n+4}$.
\end{theorem}

\begin{theorem}
\label{theorem:rnextensions}
If $\QP$ is a one-stage extension of $\R_{2^n+4}$ and $|\P| = 2$, then $\QP \cong \R_{2^{n+1}+4}$.
\end{theorem}

In this section we focus on Theorem~\ref{theorem:tnextensions}, and we prove Theorem~\ref{theorem:rnextensions} in the following section.

Throughout the discussion there will be the implicit assumption that all quotients encountered are faithful.  This is a slightly suspicious assumption, since it is unknown whether there exists an unfaithful quotient.  However, since the argument proceeds ``ground-up'' by one-stage extensions, we are safe: a careful check of the proofs reveals that every extension under consideration preserves faithfulness.  Therefore, if there exists an unfaithful quotient, it must necessarily satisfy $|\P| > 2$, and so will not interfere with the present argument.  We will not be too careful about stating and restating this assumption of faithfulness, but in all cases the checks are routine.

\subsection*{The Structure of $\T_n$}

For the remainder of this section, fix a set of games $\A$, and suppose that $\Q(\A) \cong \T_n$, where $n \geq 2$.  The structure of $\QP = \Q(\A)$ is described as follows.  $\Q = \K \cup \{1,a\}$, where $\K \cong \mathbb{Z}_2^n$ and $a^2 = 1$.  We write $\K = \{z_0,z_1,\ldots,z_{2^n-1}\}$, where $z_0$ is the identity, $z_1 = az_0$, and $z_i$ corresponds to Grundy value $i$.  

Now fix a game $G \neq 0$ with $\opts(G) \subset \A$, and write $m = \G(G)$, $\B = \cl(\A \cup \{G\})$, and $(\Q^+,\P^+) = \Q(\B)$.



\begin{definition}
Let $\E \subset \Q$.  We say that $\E$ is \emph{complemented} if $\E \cap \{a,z\} \neq \emptyset$ and $\E \cap \{1,az\} \neq \emptyset$.
\end{definition}

\begin{lemma}
\label{lemma:complementedTn}
If $\Phi''G$ is complemented, then $2n \cdot G$ is a \PPos{} for all $n \geq 1$.
\end{lemma}

\begin{proof}
Write the copies of $G$ in pairs, as $n \cdot (G+G)$.  Second player follows the mirror-image strategy on each pair \emph{until} her move would remove the last copy of $G$.  If that is the case, then the position must be
\[G + G' + Y, \textrm{ with } Y \in \A,\]
and since second player has been following the mirror-image strategy, we necessarily have $\G(Y) = 0$.

\vsp\noindent\emph{Case 1}: $\Phi(G' + Y) \in \K$.  Then second player moves to $G' + G' + Y$.  Since $\Phi(G' + G' + Y) \in \K$ and $\G(G' + G' + Y) = 0$, we necessarily have
\[\Phi(G' + G' + Y) = z \in \P.\]

\vsp\noindent\emph{Case 2}: $\Phi(G' + Y) = 1$.  Then second player chooses an $H \in \opts(G)$ with $\Phi(H) \in \{a,z\}$, as guaranteed by complementarity, and we have
\[\Phi(H + G' + Y) = \Phi(H)\cdot\Phi(G' + Y) = \Phi(H)\cdot 1 \in \P.\]

\vsp\noindent\emph{Case 3}: $\Phi(G' + Y) = a$.  Then second player chooses $H$ with $\Phi(H) \in \{1,az\}$, to the same effect.
\end{proof}

\begin{lemma}
\label{lemma:complementedTn2}
Assume that $|\P^+| \geq 2$, $\Phi''G$ is complemented, and $m = 0$ or $1$.  Fix $Y \in \A$ with $\Phi(Y) \in \K$ and $\G(Y) \neq \G(G)$.  Then $G + Y$ is an \NPos.
\end{lemma}

\begin{proof}
The $m = 0$ and $1$ cases are similar, so suppose $m = 0$.   By Lemma~\ref{lemma:complementedTn}, $G+G$ is a \PPos, so either $G+G \equiv_\B \Star$ or $G+G \equiv_\B Y+Y$.  But again by Lemma~\ref{lemma:complementedTn}, $4 \cdot G$ is a \PPos, so necsessarily $G+G \equiv_\B Y+Y$.

Now consider $G + G + G$.  A typical option is $G' + G + G \equiv_\B G' + Y + Y$; but $\G(G') \neq 0$, so
\[\Phi(G' + Y + Y) = \Phi(G')z \not\in \P.\]
Therefore $G + G + G$ is also a \PPos.  Assume (for contradiction) that $G + Y$ is also a \PPos.  Then either $G+Y \equiv_\B \Star$ or $G+Y \equiv_B Y+Y$.  But $G + G + Y + Y$ is also a \PPos, since it is equivalent to $4 \cdot Y$, so necessarily $G + Y \equiv_\B Y + Y$.  But now
\[G + G + G \equiv_\B G + Y + Y \equiv_\B Y + Y + Y,\]
a contradiction, since $Y + Y + Y$ is an \NPos.
\end{proof}

\begin{definition}
Fix $\E \subset \Q$.  The \emph{discriminant} $\Delta = \Delta(\E)$ is given by
\[\Delta = \E \cap \{1,a,z,az\}.\]
We say that $\E$ is \emph{restive} if $\Delta = \{1,z\}$ or $\{a,az\}$, \emph{restless} if $\Delta = \{a,z\}$ or $\{1,az\}$, and \emph{tame} otherwise.  We say that $\E$ is \emph{wild} if it is restive or restless.
\end{definition}

\begin{lemma}
\label{lemma:tameextensionTn}
Assume that $\Phi''G$ is tame.  If $m < 2^n$, then $\Q(\B) \cong \T_n$; if $m = 2^n$, then $\Q(\B) \cong \T_{n+1}$.  In either case, we have
\[
\Phi(G) = \begin{cases}
1 & \textrm{if $\Delta = \{a\}$}; \\
a & \textrm{if $\Delta = \{1\}$}; \\
z_m & \textrm{otherwise}.
\end{cases}\]
\end{lemma}

\begin{proof}
In each of the three cases, it is easily seen that $\Phi''G$ satisfies condition~(i) of the Generalized Mex Rule.  We now verify condition~(ii).

\vspace{0.15cm}\noindent
\emph{Case 1}: $\Delta = \{a\}$.  With $x = 1$, condition (ii) is equivalent to: for every \NPos{} $Y \in \A$, either $\Phi(Y') \in \P$ for some $Y'$, or else $x'\Phi(Y) \in \P$ for some $x' \in \E$.  But if $Y \neq 0$, then the first of these two conditions is satisfied \emph{a priori}; while if $Y = 0$, then $x' = a$ suffices for the second.

\vspace{0.15cm}\noindent
\emph{Case 2}: $\Delta = \{1\}$.  We must verify (ii) with $x = a$.  Fix $Y \in \A$ and $n \geq 0$ and suppose $a^{n+1}\Phi(Y) \not\in \P$.  If $n$ is odd, then $Y$ is an \NPos, so either $Y = 0$ or some $Y'$ is a \PPos.  If $Y = 0$, then we have $a^n \cdot 1 \cdot \Phi(Y) = a \in \P$; if $Y'$ is a \PPos, then $a^{n+1}\Phi(Y') \in \P$.  Finally, if $n$ is even, then $Y + \Star$ is an \NPos.  So either $Y$ is a \PPos, in which case $a^n \cdot 1 \cdot \Phi(Y) = \Phi(Y) \in \P$; or else $Y' + \Star$ is a \PPos, in which case $a^{n+1}\Phi(Y') = a\Phi(Y') \in \P$.

\vspace{0.15cm}\noindent
\emph{Case 3}: $\Delta \neq \{a\},\{1\}$.  Fix $Y \in \A$ and $n \geq 0$ and suppose $x^{n+1}\Phi(Y) \not\in \P$.  If $n$ is odd, then $x^{n+1} = z$, so $\Phi(Y) \neq 1,z$.  Therefore $\G(Y) \neq 0$, and $Y$ has some option $Y'$ with $\G(Y') = 0$.  Therefore $x^{n+1}\Phi(Y') = z \in \P$.

If $n$ is even, then $x^{n+1} = z_k$, so $\G(Y) \neq k$.  If $\G(Y) > k$, then there is some option $Y'$ with $\G(Y') = m$; hence $x^{n+1}\Phi(Y') = z \in \P$.  So suppose $\G(Y) < m$.  Then there is some option $G'$ of $G$ with $\G(G') = \G(Y)$.  There are three subcases.

\vspace{0.15cm}\noindent
\emph{Subcase 3a}: $n > 0$ or $\Phi(G') \in \K$ or $\Phi(Y) \in \K$.  Then we have immediately that $x^n\Phi(G')\Phi(Y) = z \in \P$.

\vspace{0.15cm}\noindent
\emph{Subcase 3b}: $n = 0$ and $\Phi(G') = \Phi(Y) = 1$.  Then $1 \in \Delta$.  Now $\Delta \neq \{1\}$ (since we are in Case 3), and furthermore $\Delta \neq \{1,az\}$ (since $\Phi''G$ is tame).  So either $a \in \Phi''G$ or $z \in \Phi''G$.  But if $x' = a$ or $z$, then $x'\Phi(Y) \in \P$, as needed.

\vspace{0.15cm}\noindent
\emph{Subcase 3c}: $n = 0$ and $\Phi(G') = \Phi(Y) = a$.  Then $a \in \Delta$.  Now $\Delta \neq \{a\}$ (since we are in Case 3), and furthermore $\Delta \neq \{a,z\}$ (since $\Phi''G$ is tame).  So either $1 \in \Phi''G$ or $az \in \Phi''G$.  But if $x' = 1$ or $az$, then $x'\Phi(Y) \in \P$, as needed.
\end{proof}

\begin{lemma}
\label{lemma:restlessextensionTn}
Assume that $\Phi''G$ is restless.  Then $|\P^+| \geq 3$.
\end{lemma}

\begin{proof}
\emph{Case 1}: $\Delta = \{1,az\}$.  Then $\{a,z\} \cap \E = \emptyset$, so $G$ is a \PPos.  Furthermore, if $G'$ is an option with $\Phi(G') = 1$ (resp.~$az$), then $\Phi(G')z \in \P$ (resp.~$\Phi(G')az \in \P$).  This shows that $G + \Star2_2$ (resp.~$G + \Star 2_3$) is an \NPos.  Therefore $G \not\equiv_\B \Star2_2$ and $G \not\equiv_\B \Star$; since $G$ is a \PPos, this implies $|\P^+| \geq 3$.

\vspace{0.15cm}\noindent
\emph{Case 2}: $\Delta = \{a,z\}$.  Then $\{1,az\} \cap \E = \emptyset$, so $\{a,z\} \cap a\E = \emptyset$, and hence $G + \Star$ is a \PPos.  Just as in Case~1, we see that $G + \Star \not\equiv_\B \Star$ and $G + \Star \not\equiv_\B \Star2_2$, so again $|\P^+| \geq 3$.
\end{proof}

\begin{lemma}
\label{lemma:restiveextensionTn}
Assume that $\Phi''G$ is restive and $|\P^+| = 2$.  Then $\Q(\B) \cong \mathcal{R}_{2^n+4}$ and
\[\Phi(G) = \begin{cases}
t & \textrm{if $\Delta = \{a,az\}$;} \\
at & \textrm{if $\Delta = \{1,z\}$.}
\end{cases}\]
\end{lemma}

\begin{proof}
The argument is similar in both cases, so suppose $\Delta = \{a,az\}$.  Now in $\mathcal{R}_{2^n+4}$ it is easy to compute $\mathcal{M}_t = \Q \setminus \{1,t,z\}$.  Since $\E \cap \{1,z\} = \emptyset$, condition~(i) of the Generalized Mex Rule is therefore trivially satisfied.

For (ii), fix $Y \in \A$ and $n \geq 0$ and suppose that $t^{n+1}\Phi(Y) \not\in \P$.  There are three cases.

\vsp\noindent\emph{Case 1}: $n > 0$.  Then $t^{n+1} = z$, so necessarily $\G(Y) > 0$.  Therefore $t^{n+1}\Phi(Y') \in \P$, where $Y'$ is any option with $\G(Y') = 0$.

\vsp\noindent\emph{Case 2}: $n = 0$ and $\Phi(Y) \not\in \K$.  If $\Phi(Y) = 1$, then we have $a\Phi(Y) \in \P$; if $\Phi(Y) = a$, then $az\Phi(Y) \in \P$.  Since $a,az \in \E$, this suffices.

\vsp\noindent\emph{Case 3}: $n = 0$ and $\Phi(Y) \in \K$.  Then $\G(Y) \neq 0$.  If $\Phi(Y') = z$ for some $Y'$, then we are done, since $t\Phi(Y') \in \P$, so assume $\Phi(Y') \neq z$ for all $Y'$.

Now since $G$ is restive, it is complemented, so by Lemma~\ref{lemma:complementedTn2} $G + Y$ is an \NPos.  Consider a typical $G + Y'$.  By assumption, $\Phi(Y') \neq z$.  If $\Phi(Y') = 1$, then $G' + Y'$ is a \PPos, where $\Phi(G') = a$.  If $\Phi(Y') = a$ or $az$, then $G' + Y'$ is a \PPos, where $\Phi(G') = az$.  If $\G(Y') \geq 2$, then by Lemma~\ref{lemma:complementedTn2} $G + Y'$ is \emph{a priori} an \NPos.  So in all cases, $G + Y'$ is an \NPos.

But $G + Y$ is an \NPos, so we must have $G' + Y$ a \PPos, for some $G'$.  Then $x'\Phi(Y) \in \P$, where $x' = \Phi(G')$, completing the proof.
\end{proof}

\begin{proof}[Proof of Theorem~\ref{theorem:tnextensions}]
Immediate from the preceding lemmas.
\end{proof}

\section{One-Stage Extensions of $\mathcal{R}_{2^n+4}$}
\label{section:Rn}

In this section we generalize much of the machinery of Section~\ref{section:Tn}.  Note that $\mathcal{R}_{2^n+4} = \mathcal{T}_n \cup \{t,at\}$, where $t^2 = tz = z$.

For the rest of this section, assume that $\QP = \Q(\A)$ is faithful, with $\QP \cong \R_{2^n+4}$.  Fix $G$ with $\opts(G) \subset \A$, and write $\B = \cl(\A \cup \{G\})$, $(\Q^+,\P^+) = \Q(\B)$, $\E = \Phi''G$, and $m = \G(G)$.

\begin{definition}
A subset $\E \subset \Q$ is said to be \emph{complemented} if $\{a,z\} \cap \E \neq \emptyset$ and $\{1,az\} \cap \E \neq \emptyset$.
\end{definition}

We can very quickly reduce to the case where $\Phi''G$ is complemented.

\begin{lemma}
\label{lemma:RnNotAT}
Assume that $\Phi''G$ is \emph{not} complemented.  If $|\P^+| = 2$, then $\Q(\B) \cong \Q(\A)$.
\end{lemma}

\begin{proof}
\emph{Case 1}: $\{a,z\} \cap \E = \emptyset$.  Since $\P = \{a,z\}$, this immediately implies that $G$ is a \PPos, so since $|\P^+| = |\P| = 2$, we must have $G \equiv_\B Y$ for some $Y \in \A$.  Therefore $\Q(\B) \cong \Q(\A)$.

\vspace{0.15cm}\noindent
\emph{Case 2}: $\{1,az\} \cap \E = \emptyset$.  If $G$ is a \PPos, then the argument is just as in Case 1.  Otherwise, consider $G + \Star$.  Since $\{1,az\} \cap \E = \emptyset$, we have $\P \cap a\E = \emptyset$, so every $G' + \Star$ is an \NPos.  Since $G + 0$ is also an \NPos, we conclude that $G + \Star$ is a \PPos.

But this implies $G + \Star \equiv_\B Y$ for some $Y \in \A$, whence $G \equiv_\B Y + \Star$, and again we have $\Q(\B) \cong \Q(\A)$.
\end{proof}

We now consider the case when $\E$ is complemented.  The key fact about complementarity is the following (cf.~Lemma~\ref{lemma:complementedTn}).

\begin{lemma}
\label{lemma:RnAlmostTame2N}
If $\Phi''G$ is complemented, then $2n \cdot G$ is a \PPos{} for all $n \geq 1$.
\end{lemma}

\begin{proof}
Identical to the proof of Lemma~\ref{lemma:complementedTn}.
\end{proof}

\begin{lemma}
\label{lemma:RnNoGY2}
Assume that $\Phi''G$ is complemented and $m = 0$ or $1$, and fix $Y \in \A$ with $\G(Y) \geq 2$.  Then $G + Y$ is an \NPos.
\end{lemma}

\begin{proof}
Identical to the proof of Lemma~\ref{lemma:complementedTn2}.
\end{proof}

\begin{lemma}
\label{lemma:RnNoGY}
Assume that $m \geq 2$, and fix $Y \in \A$ with $\Phi(Y) \in \{t,at\}$.  Then $G + Y$ is an \NPos.
\end{lemma}

\begin{proof}
First choose $G'$ with $\G(G') = \G(Y)$.  Then $\G(G'+Y) = 0$, so $\Phi(G'+Y) \in \{1,t,z\}$.  In all cases, $\Phi(G'+Y)z \in \P$, so $\Phi^+(G+Y)z^+ \not\in \P$.  Thus $\Phi^+(G+Y) \neq z^+$.

Next choose $G'$ with $\G(G') = \G(Y) \oplus 1$.  Then $\G(G'+Y) = 1$, so $\Phi(G'+Y) \in \{a,at,az\}$.  In all cases, $\Phi(G'+Y)az \in \P$, so $\Phi^+(G+Y)a^+z^+ \not\in \P$.  Thus $\Phi^+(G+Y) \neq a^+$.

Since $\P^+ = \{a^+,z^+\}$, this shows that $G + Y$ is an \NPos.
\end{proof}

\begin{lemma}
\label{lemma:RnOnlyZM}
Assume that $\Phi''G$ is complemented and $m \geq 2$, and fix $Y \in \A$.  Then $G+Y$ is a \PPos{} iff $\Phi(Y) = z_m$.
\end{lemma}

\begin{proof}
If $\Phi(Y) = z_i$, for some $i < m$, then $G' + Y$ (with $\G(G') = i$) is \emph{a priori} a \PPos, so $G + Y$ is an \NPos.  If $\Phi(Y) = z_i$ for some $i > m$, then $\G(Y') = m$ for some $Y'$, so $\Phi(Y') = z_m$.  By induction on the birthday of $Y$, we have that $G + Y'$ is a \PPos, so again $G + Y$ is an \NPos.

If $\Phi(Y) = z_m$, then by induction every $G + Y'$ is an \NPos.  Likewise, for every $G'$ we have $\G(G' + Y) \neq 0$ and $\Phi(G' + Y) \in \K$, so every $G' + Y$ is also an \NPos.  Therefore $G + Y$ is a \PPos.

This leaves only the cases $\Phi(Y) \in \{1,a,t,at\}$.  But if $\Phi(Y) = 1$ (resp.~$a$), then $G' + Y$ is a \PPos, where $\Phi(G') \in \{a,z\}$ (resp.~$\{1,az\}$), as guaranteed by complementarity.  Therefore $G + Y$ is an \NPos.  Conversely, if $\Phi(Y) \in \{t,at\}$, then Lemma~\ref{lemma:RnNoGY} guarantees that $G + Y$ is an \NPos.
\end{proof}

We now proceed with the main line of proof.  There are two fundamental cases, each stated as a separate lemma: $m \geq 2$, and $m \in \{0,1\}$.

\begin{lemma}
\label{lemma:RnAT2}
Assume that $\Phi''G$ is complemented and $m \geq 2$.  Then $\Q(\B) \cong \R_{2^n+4}$ or $\R_{2^{n+1}+4}$ and $\Phi^+(G) = z_m$.
\end{lemma}

\begin{proof}
This is much like Lemma~\ref{lemma:kernelextensions}.  It suffices to verify conditions~(i) and~(ii) in the Generalized Mex Rule.  Since $m \geq 2$, we have $\mathcal{M}_{z_m} = \Q \setminus \{z_m\}$.  Since $\G(G) = m$, this suffices for~(i).  For~(ii), fix $Y \in \A$ and $n \geq 0$, and suppose $z_m^{n+1} \Phi(Y) \not\in \P$.

If $n$ is odd, then $\G(Y) > 0$, so $z_m^{n+1} \Phi(Y') \in \P$, where $\G(Y') = 0$.

If $n$ is even, then $\G(Y) \neq m$.  If $\G(Y) > m$, then $z_m^{n+1}\Phi(Y') \in \P$, where $\G(Y') = m$.  If $2 \leq \G(Y) < m$, then let $i = \G(Y)$.  In this case $z_i$ is the unique element of $\Q$ with $\G$-value $i$, so necessarily $z_i \in \E$.  Since $z_i \cdot z_i \in \P$, this suffices.

If $\G(Y) = 0$, then we have $\Phi(Y) \in \{1,t,z\}$.  If $\Phi(Y) = z$, then $x'\Phi(Y) \in \P$ for any $x' \in \E \cap \{1,t,z\}$.  If $\Phi(Y) = 1$, then since $\E$ is complemented, we have $\E \cap \{a,z\} \neq \emptyset$; and either choice suffices.  This leaves only the case $\Phi(Y) = t$.  If $x't \in \P$ for some $x' \in \E$, then we are done.  Otherwise, $G' + Y$ is an \NPos{} for every $G'$.  But by Lemma~\ref{lemma:RnNoGY} (and the assumption $|\P^+| = 2$), we know that $G + Y$ is an \NPos, so some $G + Y'$ must be a \PPos.  By Lemma~\ref{lemma:RnOnlyZM}, we have specifically $\Phi(Y') = z_m$, whence $z_m^{n+1}\Phi(Y') = z \in \P$, as needed.

Finally, if $\G(Y) = 1$, then $\Phi(Y) \in \{a,at,az\}$, and the proof proceeds just as in the $\G(Y) = 0$ case.
\end{proof}

\begin{lemma}
Assume that $\Phi''G$ is complemented and $m = 0$ (resp.~$1$).  Then $\Q(\B) \cong \R_{2^n+4}$, and $\Phi(G) = t$ (resp.~$at$).
\end{lemma}

\begin{proof}
The two cases are essentially identical, so assume $m = 0$.  As always, we use the Generalized Mex Rule.  Note that \[\mathcal{M}_t = \Q \setminus \{1,t,z\} = \{x : \G(x) \neq 0\},\]
and since $\G(G) = 0$, this suffices for~(i).  For~(ii), fix $Y \in \A$ and $n \geq 0$, and suppose $t^{n+1}\Phi(Y) \not\in \P$.  There are four cases.

\vsp\noindent
\emph{Case 1}: $n \geq 1$.  Then $t^{n+1} = z$, so $z\Phi(Y) \not\in \P$.  Thus $\Phi(Y) \neq 1,t,z$, so $\G(Y) \neq 0$.  We conclude that $t^{n+1}\Phi(Y') \in \P$, where $Y'$ is any option with $\G(Y') = 0$.

\vsp\noindent
\emph{Case 2}: $n = 0$ and $\G(Y) = 0$.  Then $\Phi(Y) \in \{1,t,z\}$, and since $t\Phi(Y) \not\in \P$, necessarily $\Phi(Y) = 1$.  But since $\E$ is complemented, $\E \cap \{a,z\} \neq \emptyset$, so $x'\Phi(Y) \in \P$, where $x' = a$ or $z$.

\vsp\noindent
\emph{Case 3}: $n = 0$ and $\G(Y) = 1$.  Since $\E$ is complemented, we have $\E \cap \{1,az\} \neq \emptyset$.  Since $m = 0$, we know that $1 \not\in \E$, so necessarily $az \in \E$.  Since $\G(Y) = 1$, we always have $az\Phi(Y) \in \P$, so this suffices.

\vsp\noindent
\emph{Case 4}: $n = 0$ and $\G(Y) \geq 2$.  If $t \in \Phi''Y$ or $z \in \Phi''Y$, then $t\Phi(Y') = z$ and there is nothing to prove.  Otherwise, put $i = \G(Y)$; to complete the proof, it suffices to show that $z_i \in \E$, because $z_i\Phi(Y) = z \in \P$.  So consider $G + Y$.  We first show that every $G + Y'$ is an \NPos.  If $\G(Y') = 0$, then $\Phi(Y') = 1$ (since we are assuming $t,z \not\in \Phi''Y$).  Since $G$ is complemented and $\G(G) = 0$, we necessarily have $a \in \Phi''G$, so $a\Phi(Y') \in \P$ and hence $G + Y'$ is an \NPos.  If $\G(Y') = 1$, then since $G$ is complemented and $\G(G) = 0$, we necessarily have $az \in \Phi''G$, so $az\Phi(Y') \in \P$ and again $G + Y'$ is an \NPos.  Finally, if $\G(Y') \geq 2$, then the desired conclusion follows from Lemma~\ref{lemma:RnNoGY2}.

This shows that every $G + Y'$ is an \NPos.  But by Lemma~\ref{lemma:RnNoGY2}, $G + Y$ itself is an \NPos.  Therefore some $G' + Y$ is necessarily a \PPos.  Since $\Phi(Y) = z_i$, we conclude that $\Phi(G') = z_i$ as well, completing the proof.
\end{proof}

\section{Uniqueness of $\mathcal{R}_8$}
\label{section:R8}

The following theorem emerges readily from previous work.

\begin{theorem}
\label{theorem:R8unique}
$\mathcal{R}_8$ is the only mis\`ere quotient of order 8 (up to isomorphism).
\end{theorem}

\begin{proof}
Let $\QP$ be a mis\`ere quotient of order 8.  By Lemma~\ref{lemma:stages}, $\QP$ must arise as a one-stage extension of $\T_2$.  So there is some closed set $\A$, and some $G$ with $\opts(G) \subset \A$, such that
\[\Q(\A) \cong \T_2 \textrm{ and } \Q(\A \cup \{G\}) \cong \QP.\]
Let $\Phi : \cl(\A \cup \{G\}) \to \Q$ be the quotient map, and write
\[a = \Phi(\Star),\ b = \Phi(\Star 2),\ t = \Phi(G).\]
Since $\Q(\mathscr{A}) \cong \mathcal{T}_2$ and $\Q(\mathscr{A} \cup \{G\}) \not\cong \mathcal{T}_2$, $t$ is not in the submonoid generated by $a,b$.  Thus neither is $at$ (since $a^2 = 1$), and it follows immediately that
\[\Q = \T_2 \cup \{t,at\}.\]

Now put $\E = \Phi''G$.  $\E$ cannot be tame, since then Lemma~\ref{lemma:tameextensionTn} would imply that $\QP \cong \T_2$ or $\T_3$, neither of which has order 8.

If $\E$ is restive, then either $\{1,z\} \subset \E$ or $\{a,az\} \subset \E$, and it follows that $G$ and $G + \Star$ are both \NPos{}s.  Therefore $t,at \not\in \P$, so $|\P| = 2$.  By Theorem~\ref{theorem:tnextensions}, we have $\QP \cong \R_8$.

We complete the proof by assuming $\E$ to be restless and obtaining a contradiction.  There are two cases.


\vspace{0.15cm}\noindent
\emph{Case 1}: $\Delta = \{1,az\}$.  Then $a,z \not\in \E$, so $G$ is a \PPos.  Therefore $\Star 2 + G$ is an \NPos; and since $\Phi(\Star 2 + G) = bt$, we have $bt \not\in \P$, so that $bt \in \{1,b,ab,az,at\}$.  To obtain a contradiction, we show that $\Star 2 + G$ is distinguishable from some representative of each of these possibilities.

The table below summarizes.  The first column of each row lists one possibility for $bt$, along with an inequality $x \neq y$ that rules out this possibility.  In each case, $x$ is known to be in $\mathcal{P}$, and the second column exhibits an \NPos{} $Y$ that witnesses $y \not\in \mathcal{P}$.  The winning move $Y'$ is shown in the third column; the notation $\Phi^{-1}(x)$ is used to represent a typical option of $G$ with pretension $x$.

\begin{center}
\begin{tabular}{|c@{\hspace{0.25cm}$\Leftarrow$\hspace{0.25cm}}c|c|c|}
\hline
\multicolumn{2}{|c|}{Distinction(s)} & Typical \NPos & Winning Move \bigstrut \\
\hline
$1 \neq bt$ & $a \neq abt$ & & \\
$az \neq bt$ & $z \neq abt$ & $\Star + \Star 2 + G$ & $\Star + \Star + G$ \\
$at \neq bt$ & $t \neq abt$ & & \\ \hline
$b \neq bt$ & $z \neq zt$ & $\Star 2 + \Star 2 + G$ & $\Star 2 + \Star 2 + \Phi^{-1}(1)$ \bigstrut \\ \hline
$ab \neq bt$ & $z \neq azt$ & $\Star + \Star 2 + \Star 2 + G$ & $\Star + \Star 2 + \Star 2 + \Phi^{-1}(az)$ \bigstrut \\
\hline
\end{tabular}
\end{center}

\noindent
\emph{Case 2}: $\Delta = \{a,z\}$.  This is similar.  Clearly $G$ is an \NPos, so since $1,az \not\in \E$, we have that $\Star + G$ is a \PPos.  As before, this implies that $\Star 2 + G$ is an \NPos.  The following table parallels the table from Case 1.

\begin{center}
\begin{tabular}{|c@{\hspace{0.25cm}$\Leftarrow$\hspace{0.25cm}}c|c|c|}
\hline
\multicolumn{2}{|c|}{Distinction(s)} & Typical \NPos & Winning Move \bigstrut \\
\hline
$1 \neq bt$ & $a \neq abt$ & & \\
$az \neq bt$ & $z \neq abt$ & $\Star + \Star 2 + G$ & $\Star + G$ \\
$t \neq bt$ & $at \neq abt$ & & \\ \hline
$b \neq bt$ & $z \neq zt$ & $\Star 2 + \Star 2 + G$ & $\Star 2 + \Star 2 + \Phi^{-1}(z)$ \bigstrut \\ \hline
$ab \neq bt$ & $z \neq azt$ & $\Star + \Star 2 + \Star 2 + G$ & $\Star + \Star 2 + \Star 2 + \Phi^{-1}(a)$ \bigstrut \\ \hline
\end{tabular}
\end{center}

\noindent This exhausts all possibilities and completes the proof.
\end{proof}

Theorem~\ref{theorem:R8unique} can be extended: for example, $\T_3$ is the unique mis\`ere quotient of order 10.  But the proof of Theorem~\ref{theorem:R8unique} gives us pause.  The uniqueness of $\R_8$ takes shape through a somewhat subtle combinatorial analysis.  To prove the uniqueness of $\T_3$ by hand, we would need to sharpen the restless cases of Theorem~\ref{theorem:R8unique}, and then show that every one-stage extension of $\R_8$ has order $\geq 12$.  This appears to be quite a lot of work, so we now refocus our efforts on automating this sort of analysis.

\section{Valid Transition Tables}
\label{section:validity}




\emph{Transition algebras} were introduced in \cite{siegel_200Xd}, and there they proved to be useful in the study of mex functions.  We now abstract out some of their structure.

\begin{definition}
Let $\Q$ be a commutative monoid.  A \emph{transition table on $\Q$} is a subset $T \subset \Q \times \Pow(\Q)$.
\end{definition}

Note that if $\A$ is a closed set of games, then $T(\A)$ is a transition table on $\Q(\A)$.

\begin{definition}
\label{definition:valid}
Let $T$ be a transition table on a bipartite monoid $\QP$.  $T$ is said to be \emph{valid} iff the following four conditions hold.
\begin{itemize}
\item[(i)] (parity) For each $(x,\E) \in T$, we have
\[x \in \P \Longleftrightarrow \E \neq \emptyset \textrm{ and } \E \cap \P = \emptyset.\]
\item[(ii)] (completeness) For each $x \in \Q$, there is some set $\E$ such that $(x,\E) \in T$.
\item[(iii)] (closure) If $(x,\E),(y,\mathcal{F}) \in T$, then $(xy,x\mathcal{F} \cup y\E) \in T$.
\item[(iv)] (well-foundedness) There exists a map $R : \Q \to \mathbb{N}$ (a \emph{rank function} for~$\Q$) with the following property.  $R(1) = 0$, and for each $x \in \Q$, there is some $(x,\E) \in T$ such that $R(y) < R(x)$ for all $y \in \E$.
\end{itemize}
\end{definition}

We note that condition (iv) implies (ii), but nonetheless we include (ii) for clarity.  Note also that condition (iii) implies a monoid structure, so the following definition is convenient:

\begin{definition}
A transition table $T$ is a \emph{transition algebra} if it is closed (in the sense of Definition~\ref{definition:valid}(iii)).
\end{definition}

We will use the terms ``valid transition table'' and ``valid transition algebra'' interchangeably.  The main result is the following.

\begin{theorem}
\label{theorem:tt}
Let $\QP$ be a r.b.m.\ with $1 \not\in \P$.  The following are equivalent.
\begin{enumerate}
\item[(i)] There exists a closed set of games $\A$ with $\Q(\A) = \QP$;
\item[(ii)] There exists a valid transition table $T$ on $\QP$.
\end{enumerate}
\end{theorem}

\begin{proof}
(i) $\Rightarrow$ (ii): Put $T = T(\mathscr{A})$.  It is straightforward to check that $T$ is valid.  A suitable rank function is given by $R(x) = \min\{\mathrm{birthday}(G) : \Phi(G) = x\}$.

\vspace{0.15cm}\noindent
(ii) $\Rightarrow$ (i): First define, for each $x \in \Q$, a game $H_x$ as follows.  The definition is by induction on $R(x)$.  Let $(x,\E) \in T$ be such that $R(y) < R(x)$ for each $y \in \E$, and put
\[H_x = \{H_y : y \in \E\}.\]

Now define a game $H_t$ for each $t \in T$:
\[H_{(x,\E)} = \{H_y : y \in \E\}.\]

Let
\[\mathscr{A} = \cl(\{H_t : t \in T\}).\]

We claim that $\Q(\mathscr{A}) = \QP$.

Since $\QP$ is a r.b.m., it suffices (by \cite[Proposition 4.7]{siegel_200Xd}) to exhibit a surjective homomorphism $\Phi : \mathscr{A} \to \mathcal{Q}$.  Regarding $\A$ as a free commutative monoid on the generators $H_t$, we define $\Phi$ as a monoid homomorphism by
\[\Phi(H_{(x,\E)}) = x.\]
By completeness (condition (ii) in the definition of validity), $\Phi$ is surjective.  To complete the proof, we need to show that, for all $G \in \A$,
\[
\Phi(G) \in \P \Longleftrightarrow G \neq 0 \textrm{ and } \Phi(G') \not\in \P \textrm{ for any option } G'.
\]



So fix $G = H_{t_1} + \cdots + H_{t_k}$, and write $t_i = (x_i,\E_i)$.  Write $x = x_1x_2 \cdots x_k$, and denote by $x/x_i$ the product $x_1x_2\cdots x_{i-1}x_{i+1}\cdots x_k$.  Put
\[\E = \bigcup_{1 \leq i \leq k} \frac{x}{x_i}\E_i,\]
and let $t = (x,\E)$.  By closure (condition (iii) in the definition of validity), $t \in T$.  By parity (condition (i)), we have
\[x \in \P \Longleftrightarrow \E \neq \emptyset \textrm{ and } \E \cap \P = \emptyset.\]
But clearly $\Phi(G) = x$, and $\E = \Phi''G$.  This suffices except for the case when $\E = \emptyset$; but then $G$ has no options, so $\Phi(G) = 1$.  Since we assumed that $1 \not\in \P$, this completes the proof.
\end{proof}

Theorem \ref{theorem:tt} yields an algorithm for counting the number of mis\`ere quotients of order $n$: for each r.b.m.\ of order $n$, iterate over all transition tables and check whether any are valid.  This is an atrociously poor algorithm, however; even if one could effectively enumerate the r.b.m.'s of order $n$, each one admits $2^{n2^n}$ transition tables!  Theorem \ref{theorem:tt} is still important, however, since it reduces the search for mis\`ere quotients to a finite problem.

\section{Enumerating Quotients of Small Order}
\label{section:enumeration}
\suppressfloats[t]

We now show how the techniques of the previous section can be made (reasonably) efficient.  We first show that every mis\`ere quotient can be represented by a certain restricted type of transition algebra.

\begin{definition}
Let $\QP$ be a bipartite monoid.  Fix $x_1,\ldots,x_k \in \Q$, and for $0 \leq i \leq k$ let $\S_i$ be the submonoid of $\Q$ generated by $x_1,\ldots,x_i$.  We say that $x_1,\ldots,x_k$ is a \emph{construction sequence} for $\QP$ if:
\begin{enumerate}
\item[(i)] $\S_k = \Q$;
\item[(ii)] For each $i$, $x_i \not\in \S_{i-1}$;
\item[(iii)] For each $i < k$, the reduction of $(\S_i,\P \cap \S_i)$ is a mis\`ere quotient.
\end{enumerate}
\end{definition}

\begin{definition}
Let $\QP$ be a bipartite monoid.  A transition algebra $T$ on $\QP$ is said to be a \emph{minimex algebra} if there exists a construction sequence $x_1,\ldots,x_k \in \Q$ that generates $T$ in the following sense.  Write $\E_i = \mathcal{M}_{x_i} \cap \S_{i-1}$, where the $\S_i$'s are as in the previous definition.  Then $T$ is generated by
\[(x_1,\E_1),\ldots,(x_k,\E_k).\]
We say that $T$ is the minimex algebra \emph{constructed by $x_1,\ldots,x_k$}.
\end{definition}

\begin{lemma}
\label{lemma:rankedbygenerators}
Suppose $T$ is a transition algebra on a finite r.b.m.\ $\QP$.  Fix generators $x_1,\ldots,x_k \in \Q$ and suppose that, for each $i$, there is an $\E_i \subset \S_{i-1}$ such that $(x_i,\E_i) \in T$.  Then $T$ admits a rank function.
\end{lemma}

\begin{proof}
Define a map $R^* : \Q \to \mathbb{N}^k$ as follows.  For each $x \in \Q$, write
\[x = x_1^{n_1}x_2^{n_2}\cdots x_k^{n_k},\]
choosing the \emph{lexicographically least} expression on the generators $x_1,\ldots,x_k$.  Put $R^*(x) = (n_1,\ldots,n_k)$.

Now order the elements of $\mathbb{N}^k$ lexicographically.  We claim that $R^*$ is a ``rank function''  under this ordering.  For if $R^*(x) = (n_1,\ldots,n_k)$, then let
\[(x,\E) = (x_1,\E_1)^{n_1}(x_2,\E_2)^{n_2}\cdots (x_k,\E_k)^{n_k}.\]
By the assumptions on the $\E_i$, we know that $R^*(y) < R^*(x_i)$ for each $y \in \E_i$.  Therefore $R^*(y) < R^*(x)$ for each $y \in \E$.

Finally, $R^*$ can be converted into a suitable rank function $R : \mathcal{Q} \to \mathbb{N}$ by enumerating the finite range of $R^*$.
\end{proof}

\begin{theorem}
\label{theorem:minimex}
Let $\QP$ be a finitely generated r.b.m.\ with $1 \not\in \P$.  The following are equivalent.
\begin{enumerate}
\item[(i)] There exists a closed set of games $\A$ with $\Q(\A) = \QP$;
\item[(ii)] There exists a valid minimex algebra on $\QP$.
\end{enumerate}
\end{theorem}

\begin{proof}
(ii) $\Rightarrow$ (i) is immediate from Theorem~\ref{theorem:tt}, since every minimex algebra is automatically a valid transition table.  So we must prove (i) $\Rightarrow$ (ii).

Since $\Q$ is finitely generated, we may assume that $\A$ is also finitely generated (passing, if necessary, to a suitable f.g.\ subset of $\A$, and noting that the closure of a f.g.\ set is f.g.).  Choose generators $H_1,\ldots,H_l$ for $\A$ such that $\opts(H_i) \subset \<H_1,\ldots,H_{i-1}\>$ for each $i$.

Put $y_i = \Phi(H_i)$ and consider the sequence $y_1,\ldots,y_l \in \Q$.  Define a subsequence $y_{j_1},\ldots,y_{j_k}$ inductively: let $j_i$ be the \emph{least} index such that
\[y_{j_i} \not\in \S_{i-1} = \<y_{j_1},\ldots,y_{j_{i-1}}\>,\]
and stop when the subsequence $y_{j_1},\ldots,y_{j_k}$ generates $\Q$.  To avoid excessive use of nested subscripts, put $x_i = y_{j_i}$.

We claim that $x_1,\ldots,x_k$ is a construction sequence.  Conditions (i) and~(ii) are immediate from the inductive definition, and for (iii) note that
\[(\S_i,\P \cap \S_i) \textrm{ reduces to } \Q(H_1,H_2,H_3,\ldots,H_{j_i}).\]

Next let $\E_i = \Phi''H_{j_i}$ and let $U$ be the submonoid of $T(\A)$ generated by $(x_i,\E_i)$.  We claim that $U$ is valid.  Conditions (i) and (iii) (in the definition of ``valid'') are immediate, since $U$ is a submonoid of a valid transition table; and condition (ii) follows because the $x_i$'s generate $\Q$.  Finally, the choice of $x_i$'s guarantees that $\E_i \subset \S_{i-1}$, so (iv) is a consequence of Lemma~\ref{lemma:rankedbygenerators}.

Finally, let $\E_i' = \mathcal{M}_{x_i} \cap \S_{i-1}$.  Let $U'$ be generated by $(x_i,\E_i')$.  To complete the proof, we show that $U'$ is valid; then $U'$ will satisfy all the requirements of a minimex algebra.  Conditions (ii), (iii) and (iv) follow as before.  It remains to prove (i).  Now for each $i$, we know that $\E_i \subset \S_{i-1}$.  Since $U$ is valid, we have furthermore that $\E_i \subset \mathcal{M}_{x_i}$.  Therefore $\E_i \subset \E_i'$.  It follows that, whenever $(x,\E') \in U'$, then there is some $\E \subset \E'$ with $(x,\E) \in U$.

To conclude, fix any $(x,\E') \in U'$.  If $x \in \P$, then $\E \cap \P = \emptyset$ because each $\E_i' \subset \mathcal{M}_{x_i}$.  If $x \not\in \P$, then choose $\E \subset \E'$ with $(x,\E) \in U$.  Since $U$ is valid, we know that $\E \cap \P \neq \emptyset$.  Therefore $\E' \cap \P \neq \emptyset$.  This proves (i), showing that $U'$ is a minimex algebra.
\end{proof}

We now describe the algorithm for enumerating quotients of order~$n$.  Define a \emph{construction scheme} to be a tuple $(\Q,\P,x_1,\ldots,x_k)$, such that $\QP$ is a bipartite monoid and $x_1,\ldots,x_k$ is a construction sequence for~$\Q$.  A \emph{simple extension} of $(\Q,\P,x_1,\ldots,x_k)$ is a construction scheme $(\Q^+,\P^+,x_1,\ldots,x_{k+1})$ such that $\Q \subset \Q^+$ and $\P^+ \cap \Q = \P$.

It is worth emphasizing a subtle, but crucial, technicality in the definition of construction scheme.  No restrictions are placed on the b.m.\ $\QP$.  However, it is required that every \emph{proper} initial segment $(\S_i,\P \cap \S_i)$ reduce to a genuine mis\`ere quotient.  Therefore, simple extensions are meaningful only in the special case where $\QP$ is indeed a mis\`ere quotient.

By the above theorems, $\QP$ is a mis\`ere quotient if and only if there is a construction scheme $(\Q,\P,x_1,\ldots,x_k)$ such that the minimex algebra constructed by $x_1,\ldots,x_k$ is valid.  To find all mis\`ere quotients of order~$n$, we can therefore enumerate all construction schemes of order~$n$ and check which ones generate valid minimex algebras.

This method is made useful by a crucial optimization.  Built into the definition of construction sequence is the assumption that \emph{each proper initial segment reduces to a known mis\`ere quotient}.  We can therefore use the following strategy.  First, recursively compute \emph{all} mis\`ere quotients of order $< n$.  Now start with the trivial construction scheme $(\{1\},\emptyset)$.  Given a construction scheme $\Sigma = (\Q,\P,x_1,\ldots,x_k)$, consider every possible simple extension $\Sigma^+ = (\Q^+,\P^+,x_1,\ldots,x_{k+1})$ such that $|\Q^+| \leq n$.  The key is that if $|\Q^+| < n$, then $(\Q^+,\P^+)$ \emph{must} reduce to a known quotient.  If it does not, then we can discard $\Sigma^+$ from further consideration.

We have therefore reduced the search space to small simple extensions of known quotients.  Since a simple extension is just a monoid extension by a single generator, there are relatively few possibilities, and the algorithm is tractable.  It is summarized as Algorithm~\ref{algorithm:classification}.

\begin{algorithm}[tbp] 
\hrule\vspace{0.15cm}
\begin{algorithmic}[1]
\State Recursively compute all quotients of size $< n$
\State $\mathscr{X} \gets \emptyset$
\State Put the trivial construction scheme $(\{1\},\emptyset)$ into $\mathscr{X}$
\ForAll{$\Sigma = (\Q,\P,x_1,\ldots,x_k)$ \textbf{in} $\mathscr{X}$}
\State $\mathscr{Y} \gets$ the set of all simple extensions of $\Sigma$ of order $\leq n$
\ForAll{$(\Q^+,\P^+,x_1,\ldots,x_{k+1})$ \textbf{in} $\mathscr{Y}$}
\If{$|\Q^+| = n$}
\State $T \gets$ the minimex algebra on $(\Q^+,\P^+)$
\Statex \hspace{2cm} constructed by $x_1,\ldots,x_{k+1}$
\If{$(\Q^+,\P^+)$ is reduced and $T$ is valid}
\State Output $(\Q^+,\P^+)$ \Comment{It's a mis\`ere quotient}
\EndIf
\Else \Comment{$|\Q^+| \leq n-2$}
\State $(\S,\R) \gets$ the reduction of $(\Q^+,\P^+)$
\If{$(\S,\R)$ is a mis\`ere quotient}
\State Put $(\Q^+,\P^+,x_1,\ldots,x_{k+1})$ into $\mathscr{X}$
\EndIf
\EndIf
\EndFor
\EndFor
\end{algorithmic}
\vspace{0.15cm}\hrule
\caption{Classification Algorithm.}
\label{algorithm:classification}
\end{algorithm}

\newcommand{\sh}{\textrm{\raisebox{1pt}{\tiny \#}}}

\begin{figure}
\centering
\begin{tabular}{|>{$}l<{$}|>{$}r<{$}|>{$}c<{$}|>{$}c<{$}|}
\cline{2-4}
\multicolumn{1}{c|}{} & \multicolumn{1}{c|}{$\Q$} & \P & G \bigstrut \\
\hline
\S_{12} & \hangingpresent{4.4cm}{a,b,c}{a2=1,b4=b2,b2c=b3,c2=1} & \plist{1.5cm}{a,b2,ac} & \Star 2_\sh 1 \\
\hline
\S_{12}' & \hangingpresent{4.4cm}{a,b,c}{a2=1,b3=b,c2=1} & \plist{1.5cm}{a,b2,c} & \Star 2_\sh 321 \bigstrut \\
\hline
\R_{12} & \hangingpresent{4.4cm}{a,b,c,d}{a2=1,b3=b,b2c=c,c2=b2,bd=b,cd=c,d2=b2} & \plist{1.5cm}{a,b2} & \Star 2_{\sh\sh}54321 \\
\hline
& \hangingpresent{4.4cm}{a,b,c}{a2=1,b4=b2,b2c=b3,c2=b2} & \plist{1.5cm}{a,b2,c} & \Star H2_{\sh\sh} \\
\hline
& \hangingpresent{4.4cm}{a,b,c,d}{a2=1,b3=b,bc=b,c2=b2,bd=ab,d2=b2} & \plist{1.5cm}{a,b2,d} & \Star H_\sh \\
\hline
& \hangingpresent{4.4cm}{a,b,c,d}{a2=1,b4=b2,b2c=ab3,c2=abc} & \plist{1.5cm}{a,b2,c} & \Star HK2_{\sh\sh}0 \\
\hline
\end{tabular}

\vspace{0.35cm}

$H = \Star 2_{\sh\sh}321 \qquad K = \Star 2_{\sh\sh}2_\sh$
\caption{The six mis\`ere quotients of order 12.}
\end{figure}

\begin{figure}
\centering
\[\begin{array}{ccc}
\Star 2_\sh 0      & \Star (G 2_\sh)(G 2_{\sh 2} 2_\sh)   & \Star (G2_{\sh 2}2_\sh 1)(G2_{\sh 3}2_\sh 1)\\
\Star G2_\sh 32    & \Star (G 2_\sh)(G 2_{\sh 2} 2_\sh 1) & \Star 2_{\sh\sh} 4_2 54320 \\
\Star H_{\sh G}320 & \Star (G 2_\sh)(G 2_{\sh 3} 2_\sh 1) & \Star K_2K_1KG2_\sh 321\\
\end{array}\]

\vspace{0.35cm}

$G = \Star 2_\sh 320 \qquad H = \Star 2_{\sh\sh}321 \qquad K = \Star 2_{\sh\sh}2_\sh 32$
\caption{Nine games that generate non-isomorphic quotients of order 14.}
\end{figure}

\bibliography{games}

\begin{thebibliography}{1}

\bibitem{conway_1976}
J.~H. Conway.
\newblock {\em On Numbers and Games}.
\newblock A. K. Peters, Ltd., Natick, MA, second edition, 2001.

\bibitem{flammenkamp_www_octal}
A.~Flammenkamp.
\newblock {Sprague-Grundy} values of octal games.
\newblock \raggedright
  \mbox{\url{http://wwwhomes.uni-bielefeld.de/achim/octal.html}}.

\bibitem{plambeck_2005}
T.~E. Plambeck.
\newblock Taming the wild in impartial combinatorial games.
\newblock {\em INTEGERS: The Electr. J. Combin. Number Thy.}, 5(\#G05), 2005.

\bibitem{plambeck_200X}
T.~E. Plambeck.
\newblock Advances in losing.
\newblock In M.~Albert and R.~J. Nowakowski, editors, {\em Games of No
  Chance~3}, MSRI Publications. Cambridge University Press, Cambridge,
  forthcoming.
\newblock \\ \url{http://arxiv.org/abs/math.CO/0603027}.

\bibitem{siegel_200Xd}
T.~E. Plambeck and A.~N. Siegel.
\newblock Mis\`ere quotients for impartial games.
\newblock {F}orthcoming. \raggedright
  \mbox{\url{http://arxiv.org/abs/math.CO/0609825}}.

\bibitem{siegel_mqlectures}
A.~N. Siegel.
\newblock Mis\`ere {G}ames and {M}is\`ere {Q}uotients.
\newblock {L}ecture notes. \raggedright
  \mbox{\url{http://arxiv.org/abs/math.CO/0612616}}.

\end{thebibliography}
\end{document}